\documentclass{article}
\usepackage{graphicx}
\usepackage{amsmath}
\usepackage{amsfonts}
\usepackage{amsthm}
\begin{document}

\title{How to Create a New Integer Sequence}

\author{Tanya Khovanova}
\date{December 4, 2007}
\maketitle

\begin{abstract}
There are several standard procedures used to create new
sequences from a given sequence or from a given pair of sequences. In
this paper I discuss the most popular of these procedures. For each procedure, I
give a definition and provide examples based on three famous
sequences: the natural numbers, the prime numbers and the Fibonacci
numbers.  I also add my thoughts on what makes a sequence
interesting.

My goal is to help my readers invent new sequences, differentiate interesting sequences from boring ones, and better understand sequences they encounter.
\end{abstract}


\tableofcontents

\section{Introduction}

There are several standard procedures people use to create new
sequences from a given sequence or from a given pair of sequences. Most often, I use the
word ``procedure'', which is interchangeable with ``transformation'', ``operation'' or ``method''. In this
paper I discuss the most popular procedures. I am interested only in 
integer sequences, though most of the operations can be applied
to other sequences. Here is the list of my examples combined into eleven
logical groups:

\begin{enumerate}
    \item Function Acting on a Sequence Elementwise
        \begin{itemize}
            \item Adding a constant
            \item Multiplying by a constant
            \item Reverse square
            \item Delta function
        \end{itemize}
    \item Function Acting on an Index of a Sequence Element
        \begin{itemize}
            \item Adding a constant
            \item Multiplying by a constant
            \item Square
        \end{itemize}
    \item Composition of Sequences
        \begin{itemize}
            \item Self-composition
            \item Composition
        \end{itemize}
    \item Compositional Inverse
        \begin{itemize}
            \item Left inverse
            \item Right inverse
        \end{itemize}
    \item Function Acting on Two Sequences Elementwise
        \begin{itemize}
            \item Sum of two sequences
            \item Square
            \item Product of two different sequences
            \item Concatenation of a sequence element with its reverse
            \item Concatenation of one sequence element with the reverse of another sequence element
            \item Function acting on many sequences elementwise
        \end{itemize}
    \item Set Operations
        \begin{itemize}
            \item Complement
            \item Intersection
            \item Union
        \end{itemize}
    \item Function Acting on Sets
        \begin{itemize}
            \item Reverse square
            \item Sum of two sets
            \item Product of two sets
        \end{itemize}
    \item Discrete Calculus
        \begin{itemize}
            \item Partial sums
            \item First difference
            \item Partial products
        \end{itemize}
    \item Geometric Inverse Sequence
        \begin{itemize}
            \item Geometric inverse
            \item Indicator
            \item Reverse indicator
        \end{itemize}
    \item Convolution of Two Sequences
        \begin{itemize}
            \item Self-convolution
            \item Convolution
            \item Convolutional inverse
        \end{itemize}
    \item Binomial Transform
        \begin{itemize}
            \item Binomial transform
            \item Inverse binomial transform
            \item Binomial transform III
        \end{itemize}
\end{enumerate}

For every operation I provide examples using three basic sequences. My ``lab rats'' are three very popular sequences: the natural numbers, the prime numbers and the Fibonacci numbers:

\begin{itemize}
\item 1, 2, 3, 4, 5, 6, 7, 8, 9, 10, 11, $\dots$ \newline
A000027: The natural numbers.
\item 2, 3, 5, 7, 11, 13, 17, 19, 23, 29, 31, 37, 41, $\dots$ \newline
A000040: The prime numbers.
\item 1, 1, 2, 3, 5, 8, 13, 21, 34, 55, 89, 144, 233, $\dots$ \newline
A000045: Fibonacci numbers.
\end{itemize}

For every new sequence produced I check whether this sequence is in the Online
Encyclopedia of Integer Sequences (OEIS) \cite{OEIS}. If it is there, I
provide the corresponding link and the definition from the OEIS.

After
showing the examples I discuss them. One of the topics of the
discussion is the interestingness of the results: are the produced
sequences interesting or not, and why or why not. As a regular submitter of sequences to the OEIS database, I've built an emotional isomorphism of a sequence being interesting and a sequence being worth submitting to the database. In this paper I use notions of being interesting and worth submitting interchangeably.

It is easy to create an
infinite number of sequences by combining or repeating the discussed procedures.
Most such sequences will not be interesting, and as such will not
deserve to be submitted to the OEIS, never mind
the time commitment required to submit an infinite
number of sequences. It is not very beneficial just to take a random
sequence, apply a random procedure discussed here and submit the
resulting sequence. It might be even less beneficial to take all the
sequences in the database, apply all these procedures and submit all the
results. It is very beneficial for the software that is used in the
OEIS database to incorporate these procedures while looking up a
sequence. A program of this kind exists and it is called Superseeker. Superseeker
uses its own list
of transformations which partially overlaps with the list of
procedures I discuss here. I hope in the future Superseeker will
become more powerful and will include more procedures from this list.

Meanwhile, if you take a random sequence and apply a random procedure,
there are many things that might make your new sequence very
interesting and worth submitting:

\begin{itemize}
\item Your new sequence amazes you
\item You find an extra property for your new sequence
\item You can prove something non-trivial about your new sequence.
\end{itemize}

The procedures I discuss are interesting not only for generating new
sequences, but also for decomposing existing
sequences into simpler sequences. At the end of this paper I give two
examples of building advanced sequences from my "lab rat" sequences
using the discussed procedures.

\textbf{Formalities.} I denote my main sequence to play with as
$a(n)$, where $n$ is the index. I assume that the index
$n$ starts with 1. This paper was synchronized with the OEIS in July
2007. Sequences in the OEIS might start with a different index.

\section{Function Acting on a Sequence Elementwise}

Suppose $f$ is a function from integers to integers. Then, given
a sequence $a$, we can define a sequence $b$: $b =
f(a)$; where, for each index $n$, $b(n) = f(a(n))$. That
is, each element of $b$ is equal to the function $f$ applied
to the same-indexed element of $a$. We say that the sequence
$b$ is the function $f$ acting on the sequence $a$.

There are two special cases to keep in mind. First case: if $f(n) =
n$, then $f(a(n)) = a(n)$, that is $f(n)$ acts as the
identity. Second case: if $f(n) = c$, then $b(n) = f(a(n)) =
c$. For my examples I consider four different cases for the function
$f$. The first two cases are the most standard ones: adding a
constant and multiplying by a constant. The third case is a more
complicated function I have chosen at random. The fourth case is a
delta function, which plays a special role in this paper.

\textbf{Adding a constant:} Let $f(k) = k+m$, where $m$ is an
integer. Then $b(n) = a(n) + m$. For my example sequences I
consider the special case $m = 1$:

\begin{itemize}
    \item If $a(n)$ are the natural numbers, then $b(n)$ is: \newline
    2, 3, 4, 5, 6, 7, 8, 9, 10, 11, 12, $\dots$ \newline
    A020725: Integers $\ge$ 2.
    \item If $a(n)$ are the prime numbers, then $b(n)$ is: \newline
    3, 4, 6, 8, 12, 14, 18, 20, 24, 30, 32, 38, 42, $\dots$ \newline
    A008864: Primes + 1.
    \item If $a(n)$ are the Fibonacci numbers, then $b(n)$ is: \newline
    2, 2, 3, 4, 6, 9, 14, 22, 35, 56, 90, 145, 234, $\dots$ \newline
    A001611: Fibonacci numbers (A000045) + 1.
\end{itemize}

\textbf{Discussion. Natural numbers.} Adding a constant to the natural
numbers produces a shift --- the new sequence is essentially the
same sequence as the natural number sequence itself.

\textbf{Discussion. Prime numbers.} Adding constants to the prime
numbers generates an infinite number of sequences. Which of the
constants are more interesting to add? In general, people find that
adding a very small number, like 1 or 2, is more interesting than
adding a random big number, like 117. Sometimes a
particular number exists, related to the sequence, which is especially suitable for
addition. In case of prime numbers, the number 2 is the difference
between pairs of twin primes. Therefore, I find adding the
number 2 to the prime numbers more interesting than adding the number
2 to a random sequence.

\textbf{Discussion. Fibonacci numbers.} The Fibonacci sequence is a linear
recurrence of the second order. A constant sequence is a linear
recurrence of the first order. Hence, we can expect that adding a
constant to the Fibonacci sequence can create something interesting
too. It is easy to see that if $b(n) = Fibonacci(n) + c$, then
$b(n) = b(n-1) + b(n-2) - c = 2b(n-1) - b(n-3)$. That
is, $b(n)$ is a linear recursive sequence of the third order. As a
result, I find adding a random constant to the Fibonacci sequence to
be more interesting than adding the same constant to the sequence of
prime numbers.

\textbf{Multiplying by a constant:} Let $f(k) = mk$, where
$m$ is an integer. Then $b(n) = ma(n)$. For my examples I
consider the special case $m = 2$:

\begin{itemize}
    \item If $a(n)$ are the natural numbers, then $b(n)$ is: \newline
    2, 4, 6, 8, 10, 12, 14, 16, 18, 20, 22, 24, 26, 28, 30, $\dots$ \newline
    A005843: The even numbers: a(n) = 2n.
    \item If $a(n)$ are the prime numbers, then $b(n)$ is: \newline
    4, 6, 10, 14, 22, 26, 34, 38, 46, 58, 62, 74, 82, $\dots$ \newline
    A100484: Even semiprimes.
    \item If $a(n)$ are the Fibonacci numbers, then $b(n)$ is: \newline
    2, 2, 4, 6, 10, 16, 26, 42, 68, 110, 178, 288, 466, $\dots$ \newline
    Almost A006355: Number of binary vectors of length n containing no singletons.\newline
    Also almost A055389: a(0)=1, then twice the Fibonacci sequence.\newline
    Also almost A118658: L\_n - F\_n where L\_n is the Lucas Number and F\_n is the Fibonacci Number.
\end{itemize}

\textbf{Discussion. Prime numbers.} The prime number sequence is a set
of numbers that share a special multiplicative property --- they
are all divisible only by 1 and the number itself. Because of that,
multiplication by a number might be more interesting than adding a
number to this sequence. If you look at the results for the prime
numbers, you will see that multiplication by 2 gives a new sequence
with its own description: even numbers that are a product of two
primes. At the same time adding one to the prime numbers gives a sequence
that is described in terms of this exact operation: primes plus
1. That is, I find that multiplying the prime numbers by a constant is
more interesting in general than adding a constant to the prime
numbers.

\textbf{Discussion. Fibonacci numbers.} Now let's look at the Fibonacci
sequence, which is very different from the sequence of prime
numbers. In particular, the Fibonacci sequence is a linear recursive
sequence of the second order. Because of that, when the Fibonacci sequence
is multiplied by a number the recurrence relation is preserved. The
resulting sequence keeps most of the properties of the Fibonacci
sequence. In some sense, the new sequence is almost as interesting as
the Fibonacci sequence. At the same time, there is not much need to
study two different sequences separately when they have the same recurrence relation. It is enough to study one of them, and then transfer the
properties to the other. For historical reasons the Fibonacci sequence
is the sequence of choice to study the recurrence $a(n) = a(n-1) +
a(n-2)$.

\textbf{Reverse square:} The function $f(k)$ could be any obscure function. Let $f(k) = Reverse(k^2)$. For example:

\begin{itemize}
    \item If $a(n)$ are the natural numbers, then $b(n)$ is: \newline
    1, 4, 9, 61, 52, 63, 94, 46, 18, 1, 121, 441, 961, 691, $\dots$ \newline
    A002942: Squares written backwards.
    \item If $a(n)$ are the prime numbers, then $b(n)$ is: \newline
    4, 9, 52, 94, 121, 961, 982, 163, 925, 148, 169, 9631, $\dots$ \newline
    Almost A060998: Squares of 1 and primes, written backwards.
    \item If $a(n)$ are the Fibonacci numbers, then $b(n)$ is: \newline
     1, 1, 4, 9, 52, 46, 961, 144, 6511, 5203, $\dots$ \newline
     This sequence is not in the database. It is not interesting enough to be in the database. Actually, the most interesting thing about this sequence might be its presence on this page.
\end{itemize}

\textbf{Discussion.} The result of reversing a number depends on the
base in which the number is written. We use base ten mostly because we
have 10 fingers on our hands. If we had 14 fingers, the reverse
operation would have had a very different result. For this reason,
many mathematicians feel that the reversing operation is not a
mathematical operation, and shouldn't be considered interesting or
worth looking at. This is the same reason why non-base related
submissions to the OEIS database are more encouraged than base
related submissions. At the same time, there are many sequences in the
database that are base-related and people continue submitting
them. One thing in favor of such sequences is that they often have
very short and simple descriptions. I find sequences that can be
described in two words very appealing. Also, people love symbols and
wonder about symbolic properties of numbers. We can say that
base-related sequences reflect not just properties of numbers, but also
properties of the symbols representing them.

\textbf{Delta function:} A very special case for $f$ is a
delta function. Namely $\delta_m(n) = 1$, if $n =
m$ and 0 otherwise. Let us consider an example where $m = 1$:

\begin{itemize}
    \item If $a(n)$ are the natural numbers, then $b(n)$ is: \newline
    1, 0, 0, 0, 0, 0, 0, 0, 0, 0, 0, 0, 0, 0, 0, 0, 0, 0, 0, $\dots$ \newline
    A063524: Characteristic function of 1.
    \item If $a(n)$ are the prime numbers, then $b(n)$ is: \newline
    0, 0, 0, 0, 0, 0, 0, 0, 0, 0, 0, 0, 0, 0, 0, 0, 0, 0, $\dots$ \newline
    A000004: The zero sequence.
    \item If $a(n)$ are the Fibonacci numbers, then $b(n)$ is: \newline
    1, 1, 0, 0, 0, 0, 0, 0, 0, 0, 0, 0, 0, 0, 0, 0, 0, 0, $\dots$ \newline
    A019590: Fermat's Last Theorem: a(n) = 1 if x\^n+y\^n=z\^n has a nontrivial solution in integers, otherwise a(n) = 0.
\end{itemize}

\textbf{Discussion.} Obviously, $\delta_m(n)$ acting on
a sequence $a(n)$ is equal to 1, for $n$ such that $a(n) =
m$, and is equal to 0 otherwise. If our sequence $a(n)$ never
reaches the value $m$, then the resulting sequence is the zero
sequence. The procedure of a delta function acting on a sequence can be
especially interesting if our initial sequence reaches the value
$m$ an infinite number of times.

\section{Function Acting on an Index of a Sequence Element}

Suppose $f$ is a function from integers to integers. Suppose
further that $f(n)$ is positive for positive $n$, so that
$f(n)$ is a valid index. Then, given a sequence $a$, we can
define a sequence $b$: $b = a(f)$; where, for each index
$n$, $b(n) = a(f(n))$. We say that the sequence $b$ is
the function $f$ acting on the index of the sequence $a$.

There are two special cases to keep in mind. First case: if $f(n) =
n$, then $a(f(n)) = a(n)$, that is $f(n)$ acts as the
identity. Second case: if $f(n) = c$, then $b(n) = a(f(n)) =
a(c)$. For my examples I consider 3 different cases for the
function $f$. The first two cases are the most standard ones:
adding a constant and multiplying by a constant. The third case is a
more complicated function I have chosen at random.

\textbf{Adding a constant:} Let $f(k) = k+m$, where $m$ is an
integer. Then $b(n) = a(n+m)$; that is, $b(n)$ is the same
sequence as $a(n)$ but shifted by $m$. Or, in other words,
the same sequence starting from a different place.

\textbf{Multiplying by a constant:} Let $f(k) = mk$, where $m$
is an integer. Then $b(n) = a(mn)$. For my examples I consider
the special case $m = 2$:

\begin{itemize}
    \item If $a(n)$ are the natural numbers, then $b(n)$ is: \newline
    2, 4, 6, 8, 10, 12, 14, 16, 18, 20, 22, 24, 26, 28, 30, $\dots$ \newline
    A005843:  The even numbers.
    \item If $a(n)$ are the prime numbers, then $b(n)$ is: \newline
    3, 7, 13, 19, 29, 37, 43, 53, 61, 71, 79, 89, 101, $\dots$ \newline
    A031215:  (2n)-th prime.
    \item If $a(n)$ are the Fibonacci numbers, then $b(n)$ is: \newline
    1, 3, 8, 21, 55, 144, 377, 987, 2584, 6765, 17711, $\dots$ \newline
    A001906:  F(2n) = bisection of Fibonacci sequence: a(n)=3a(n-1)-a(n-2). 
\end{itemize}

\textbf{Discussion.} You can notice that the bisection of the Fibonacci
sequence has a recurrence relation in its own right. This is not a
coincidence. In fact, if a sequence $a(n)$ satisfies the
recurrence relation $a(n) = pa(n-1) + qa(n-2)$, then its
bisection $b(n) = a(2n)$ satisfies the recurrence relation
$b(n) = (p^2 + 2q)b(n-1) - q^2 b(n-2)$.

\textbf{Square:} $f(k)$ could be any function. For example, let $f(k) = k^2$:

\begin{itemize}
    \item If $a(n)$ are the natural numbers, then $b(n)$ is: \newline
    1, 4, 9, 16, 25, 36, 49, 64, 81, 100, 121, 144, 169, $\dots$ \newline
    A000290:  The squares. 
    \item If $a(n)$ are the prime numbers, then $b(n)$ is: \newline
    2, 7, 23, 53, 97, 151, 227, 311, 419, 541, 661, 827, 1009, $\dots$ \newline
    A011757:  prime\_(n\^2).
    \item If $a(n)$ are the Fibonacci numbers, then $b(n)$ is: \newline
    1, 3, 34, 987, 75025, 14930352, 7778742049, 10610209857723, $\dots$ \newline
    A054783:  (n\^2)-th Fibonacci number.
\end{itemize}

\textbf{Discussion.} The choice of the function acting on the index in
this case is not related to the inner properties of the prime numbers
or the Fibonacci numbers. That is why many mathematicians might find
the sequences A011757 and A054783
not very interesting. Indeed, if you look at these sequences in the
database you will see that though they were submitted a long time ago,
they've received no comments. Still, I know three things that could give
extra points to the interestingness score of these two sequences:

\begin{enumerate}
    \item The sequences have short descriptions, which is always a plus.
    \item The sequences are increasing, which means they are not at random and have some order.
    \item The growth rate for these two sequences can be easily
approximated. Indeed, if $g(n)$ describes the growth rate of
$a(n)$, then $g(f(n))$ describes the growth rate of
$a(f(n))$; and it is very easy to plug the square function
into the growth rates for the prime numbers and the Fibonacci
numbers.
\end{enumerate}

\section{Composition of Sequences}

A sequence can be viewed as a function from positive integers to
integers. Vice versa, any function on integers, when restricted to the
positive integers, forms a sequence. Suppose we have two functions
$f(n)$ and $g(n)$. The function $h(n) = f(g(n))$ is
called the composition of the two functions $f$ and $g$. The
idea of function composition can be expanded to sequences. Suppose we
have two sequences: $a(n)$ and $b(n)$. Additionally, suppose
$b(n)$ is positive for every $n$. Then the sequence $c(n)
= a(b(n))$ is called the composition sequence of $a$ and
$b$.

As we've seen before, composition with the natural numbers doesn't change the sequence. That is, the natural number sequence acts as the identity for this operation.

\textbf{Note.} The composition of sequences procedure is very
similar to the two previous procedures: function acting on a sequence
element and function acting on an index. To look at this similarity in
more detail, let us start with two sequences $a(n)$ and
$b(n)$ and correspond to them two functions on positive integers:
$f(n)$ and $g(n)$. Suppose the sequence $b(n)$ is
positive, then $g(n)$ is positive. Now the composition sequence
$c(n) = a(b(n))$ is the same sequence as the function $f$
acting on elements of $b$ and also the same sequence as the
function $g$ acting on indices of $a$. On the other hand, if
the sequence $b(n)$ is not positive, we still can have a function
acting on it. In this sense, a function acting on a sequence is a more
general operation than the composition of two sequences. But for
positive sequences, both a function acting on a sequence element and a
function acting on an index are the same procedure as the composition
of two sequences.

\textbf{Self-composition:} For this special case, let $a(n) = b(n)$, then $c(n) = a(a(n))$:

\begin{itemize}
    \item If $a(n)$ are the natural numbers, then $c(n)$ is: \newline
    1, 2, 3, 4, 5, 6, 7, 8, 9, 10, 11, $\dots$ \newline
    A000027:  The natural numbers.
    \item If $a(n)$ are the prime numbers, then $c(n)$ is: \newline
    3, 5, 11, 17, 31, 41, 59, 67, 83, 109, 127, 157, $\dots$ \newline
    A006450:  Primes with prime subscripts. 
    \item If $a(n)$ are the Fibonacci numbers, then $c(n)$ is: \newline
    1, 1, 1, 2, 5, 21, 233, 10946, 5702887, 139583862445, $\dots$ \newline
    A007570:  F(F(n)), where F is a Fibonacci number.
\end{itemize}

\textbf{Discussion.} Similar to the discussion in the previous chapter,
we can give extra interestingness points to the sequences A006450
and A007570:
for their short descriptions, for being increasing sequences, for
the ease of calculating their growth rates. Is there anything else?
One can hope that if two sequences are related to each other their
composition might be an exciting sequence. A sequence is definitely
related to itself --- is this enough? Obviously, the
self-composition can't be equally interesting for every sequence. What
kinds of sequences allow the self-composition to produce something
special? Are the prime numbers and the Fibonacci numbers the best
choices to plug into the self-composition? I am not sure. I might
prefer to plug the square sequence into the self-composition:

\begin{itemize}
    \item If $a(n)$ are the squares, then $c(n)$ is: \newline
    1, 16, 81, 256, 625, 1296, 2401, 4096, 6561, 10000, $\dots$ \newline
    A000583: Fourth powers.
\end{itemize}

\textbf{Composition:} Due to triviality, I omit the cases where one of the sequences is the natural number sequence:

\begin{itemize}
    \item If $a(n)$ are the prime numbers and $b(n)$ are the Fibonacci numbers, then $c(n)$ is: \newline
    2, 2, 3, 5, 11, 19, 41, 73, 139, 257, 461, 827, 1471, $\dots$ \newline
    A030427:  Prime(Fibonacci(n)).
    \item If $a(n)$ are the Fibonacci numbers and $b(n)$ are the prime numbers, then $c(n)$ is: \newline
    1, 2, 5, 13, 89, 233, 1597, 4181, 28657, 514229, 1346269, $\dots$ \newline
    A030426:  Fibonacci(prime(n)).
\end{itemize}

\textbf{Discussion.} The prime and the Fibonacci numbers are seemingly
unrelated to each other. As a result, the following fact becomes
amazing: every Fibonacci number $F(n)$ that is prime has a prime
index $n$, with the exception of $F(4) = 3$. That means the
sequence A030426
above contains all prime Fibonacci numbers except 3. I find this
sequence very interesting.

\section{Compositional Inverse}

As I mentioned before, the sequence of natural numbers acts as the
identity under the composition operation. When we have an operation
with an identity we usually try to define an inverse object. For many mathematical operations the inverse is unique, or in the worst case there are two inverses: left and right. With sequences everything is worse than the worst case. We will see that the inverses are not always defined and there can be many of them. Let us try to bring some order to this chaos of compositional inverses of sequences.

We can start with standard definitions for left and right inverses. Namely,
given a sequence $a(n)$ we say that a sequence $b(n)$ is a
left inverse of $a$ if the sequence $b(a(n))$ is the natural
number sequence. I denote a left inverse sequence as
$leftInv(n)$. Correspondingly, a right inverse sequence is
denoted as $rightInv(n)$, and it satisfies the property that the
composition sequence $a(rightInv(n))$ is the natural numbers
sequence. It goes without saying that the sequences $leftInv(n)$ 
and $rightInv(n)$ depend on the sequence $a$. I will sometimes use the notation $leftInv(a)(n)$ 
and $rightInv(a)(n)$ in cases where I need this dependency to be explicit.

\textbf{Left inverse.} First we assume that $a(n)$ is
positive. Next, if $a(n)$ takes the same value for two different
indices $n$, then the left inverse sequence cannot be
defined. If $a(n)$ doesn't reach a number $K$ for any index
$n$, then $leftInv(K)$ could be any number. That is, in this
case the left inverse isn't defined uniquely. From here, we see that
we can define the left inverse uniquely only if $a(n)$ is a
permutation of natural numbers and in this case the left inverse
sequence is the reverse permutation.

Many of the interesting sequences are increasing. To be able to define a left inverse for an increasing sequence, we need this sequence not to take the same value for different indices. This requirement translates into a simple condition: our increasing sequence has to be strictly increasing. Suppose $a(n)$ is a strictly increasing sequence. In this case the left inverse sequence can be defined. It still might
not be unique, or more precisely, it is guaranteed not to be unique unless $a$ is the sequence of natural numbers.

Each time the left inverse is not unique we have infinitely many left inverses. To enjoy some order in this chaos of left inverses I would like to restrict candidates for the left inverse to
non-decreasing sequences. In
this case we can define two special left inverse sequences:
$minimalLeftInv$ and $maximalLeftInv$, called the minimal
left inverse and the maximal left inverse correspondingly. We define
them so that for any non-decreasing sequence $b(n)$, such that
$b(n)$ is a left inverse of $a(n)$, the following equations
are true: $minimalLeftInv(n) \le b(n) \le
maximalLeftInv(n)$.

It is easy to see that the
$minimalLeftInv(a)(n)$ is the number of elements in $a(n)$ that
are less than or equal to $n$. Also the $maximalLeftInv(a)(n)$
is the number of elements in $a(n)$ that are less than $n$,
plus 1. In particular, $maximalLeftInv(n) -
minimalLeftInv(n)$ equals 0 if $n$ belongs to $a(n)$ and
1 otherwise. From here trivially we get the following equations:
$minimal\-Left\-Inv(n) + 1 = 
maximalLeftInv(n) + \text{ characteristic function of } a(n) =
maximalLeftInv(n+1)$.

\textbf{Left inverse:}

\begin{itemize}
    \item If $a(n)$ are the natural numbers, then $leftInv(n)$ is uniquely defined and is: \newline
    1, 2, 3, 4, 5, 6, 7, 8, 9, 10, 11, $\dots$ \newline
    A000027:  The natural numbers.
    \item If $a(n)$ are the prime numbers, then $leftInv(n)$
    is any sequence such that if $n$ is a prime number then
    $leftInv(n)$ is $\pi(n)$ --- the number of primes
    less than or equal to $n$. Here are the minimal left inverse
    and the maximal left inverse:
    \begin{itemize}
        \item 0, 1, 2, 2, 3, 3, 4, 4, 4, 4, 5, 5, 6, 6, 6, 6, 7, 7, 8, 8, $\dots$ \newline
        A000720: pi(n), the number of primes $\le$ n.
        \item 1, 1, 2, 3, 3, 4, 4, 5, 5, 5, 5, 6, 6, 7, 7, 7, 7, 8, 8, 9, $\dots$ \newline
        Almost A036234: Number of primes $\le$ n, if 1 is counted as a prime.
    \end{itemize}
    
    \item If $a(n)$ are the Fibonacci numbers, then
    $leftInv(n)$ can't be defined because $a(1) = a(2) =
    1$. Suppose we start the Fibonacci sequence from the second
    1. For this slightly trimmed Fibonacci sequence we can define the
    minimal left inverse and the maximal left inverse:
    \begin{itemize}
        \item 1, 2, 3, 3, 4, 4, 4, 5, 5, 5, 5, 5, 6, 6, 6, 6, 6, 6, 6, 6, 7, $\dots$ \newline
        A072649: n occurs A000045(n) times.
        \item 1, 2, 3, 4, 4, 5, 5, 5, 6, 6, 6, 6, 6, 7, 7, 7, 7, 7, 7, 7, 7, 8, $\dots$ \newline
        A131234: Starts with 1, then n appears Fibonacci(n-1) times. 
    \end{itemize}    
\end{itemize}

\textbf{Discussion.} If you compare the descriptions above of the
minimal/maximal left inverses for primes with the minimal/maximal left
inverses for the Fibonacci sequence, you can notice a discrepancy in
these descriptions. To explain this discrepancy, let me give you other
definitions of the minimal/maximal left inverse sequences. Namely,
given an increasing sequence $a(n)$, the minimal left inverse
sequence can be described as: $n$ appears $a(n+1) -
a(n)$ times. Correspondingly, the maximal left inverse sequence can
be described as: $n$ appears $a(n) - a(n-1)$
times. With these definitions the discrepancy is explained by the fact
that for the Fibonacci sequence the expressions $a(n+1) -
a(n)$ and $a(n) - a(n-1)$ can be simplified into
$a(n-1)$ and $a(n-2)$ correspondingly.

\textbf{Right inverse.} Again we assume that $a(n)$ is
positive. It is easy to see that if $a(n)$ doesn't reach a number
$K$ for any index $n$, then the right inverse can't be
defined. If $a(n)$ takes the same value for two or more different
indices $n$, then the right inverse sequence can reach only one
of those index values (and we can choose which one). From here, we see
that we can define the right inverse uniquely only if $a(n)$ is a
permutation of natural numbers and in this case the right inverse
sequence is the reverse permutation.

Suppose $a(n)$ is a sequence that reaches every natural number
value. Therefore, the right inverse sequence can be defined. The right
inverse sequence might not be unique, but we can try to define two
special right inverse sequences: $minimalRightInv$ and
$maximalRightInv$, called the minimal and the maximal right
inverse correspondingly. We define them so that for any sequence $b(n)$, such that $b(n)$ is a right inverse of
$a(n)$, the following equations are true: $minimalRightInv(n)
\le b(n) \le maximalRightInv(n)$. It is easy to see that
$minimalRightInv(n)$ is the smallest index $k$, such that
$a(k) = n$. Also, $maximalRightInv(n)$ is the largest index
$k$, such that $a(k) = n$. It is easy to see that the
minimal right inverse is always defined. At the same time, for the
maximal right inverse to be defined, it is necessary and sufficient
that $a(n)$ reaches every value a finite number of times.

Suppose that $a(n)$ is a non-decreasing sequence that reaches
every natural number value a finite number of times. Then the maximal
right inverse is defined and

$$maximalRightInv(n) =
minimalRightInv(n+1)-1.$$

Suppose $a(n)$ is a strictly increasing sequence. Then both the minimal and maximal left inverses are defined. Moreover, both of them are non-decreasing sequences that reach every value a finite number of times. This means that we can define the minimal and maximal right inverses on the sequences $minimal\-Left\-Inv(a)(n)$ and $maximalLeftInv(a)(n)$. The following properties are true:

\begin{itemize}
     \item $minimalRightInv(minimalLeftInv(a))(n) = a(n)$
     \item $maximalRightInv(minimalLeftInv(a))(n) + 1 = a(n+1)$
     \item $minimalRightInv(maximalLeftInv(a))(n+1) = a(n) + 1$
     \item $maximalRightInv(maximalLeftInv(a))(n) = a(n)$
\end{itemize}

\textbf{Examples.} The right inverse sequence for the natural
numbers is uniquely defined and is the sequence of natural
numbers. The right inverse sequence for the prime or the Fibonacci
numbers cannot be defined. Not to leave you without an example, let
us see what happens if we make the composition of $leftInv(\textrm{prime numbers})$ with the Fibonacci numbers:

\begin{itemize}
    \item The composition of $leftInv(\textrm{prime numbers})$ with the Fibonacci numbers:\newline
    0, 0, 1, 2, 3, 4, 6, 8, 11, 16, 24, 34, 51, 74, 111, 166, 251, 376, $\dots$ \newline
    A054782:  Number of primes $\le$ n-th Fibonacci number.
\end{itemize}

\section{Function Acting on Two Sequences Elementwise}

Suppose $f = f(x,y)$ is an integer function of two integer
variables. Then, given two sequences $a(n)$ and $b(n)$, we
can define a sequence $c: c = f(a,b)$, where, for each index
$n$, $c(n) = f(a(n),b(n))$. We say that the sequence
$c$ is the function $f$ acting on the sequences $a$ and
$b$ elementwise.

This section is a generalization of the section ``Function Acting on a Sequence Elementwise''. Following my pattern from that section I consider 3 different cases for the function
$f$: the sum of two sequences, the product of two sequences and a
random function. At the same time I am breaking the pattern of the previous chapters: for the first time I am discussing how to create a new sequence using a pair of known sequences. This is the time to create a new pattern. The new pattern is the following: each time I create a new sequence based on a pair of sequences $a$ and $b$ I will look separately at two subcases. The first subcase is when $a$
is the same as $b$ and the second subcase is when they are different.

If the two sequences are the same: $a = b$, then $f(a(n),b(n))$ becomes a
function of $a(n)$. Therefore, this subcase is a particular case of
a function acting on a sequence. You might think that I have a right to skip this subcase as it formally belongs to another section of this paper. I am dropping this right in favor of fun, so this subcase stays.

\textbf{Sum of two sequences.} The sum, $s(n)$, of two sequences
$a(n)$ and $b(n)$ is defined as $s(n) = a(n) +
b(n)$. Summing a sequence with itself is the same as multiplying
this sequence by 2. We already discussed this example before, hence, we can proceed with examples
of the sums of two different basic sequences:

\begin{itemize}
    \item If $a(n)$ are the natural numbers and $b(n)$ are the prime numbers, then $s(n)$ is: \newline
    3, 5, 8, 11, 16, 19, 24, 27, 32, 39, 42, 49, 54, 57, 62, $\dots$ \newline
    A014688:  a(n) = n-th prime + n.
    \item If $a(n)$ are the natural numbers and $b(n)$ are the Fibonacci numbers, then $s(n)$ is: \newline
    2, 3, 5, 7, 10, 14, 20, 29, 43, 65, 100, 156, 246, 391, 625, $\dots$ \newline
    A002062:  n-th Fibonacci number + n.
    \item If $a(n)$ are the prime numbers and $b(n)$ are the Fibonacci numbers, then $s(n)$ is: \newline
    3, 4, 7, 10, 16, 21, 30, 40, 57, 84, 120, 181, 274, 420, 657, $\dots$ \newline
    A004397:  n-th prime + n-th Fibonacci number.
\end{itemize}

\textbf{Discussion.} I would like to introduce the not very precise idea
of a shiftable sequence. I call a sequence shiftable if it keeps some
of its properties when started from a different index. In
particular, it means that the order in which the sequence is presented
is important and is related to the sequence's properties. I consider
the prime number sequence not very shiftable: the prime numbers do not
relate to each other very well. The Fibonacci sequence is very
shiftable. If you start the Fibonacci numbers from any place in the
Fibonacci sequence, you will get a sequence with the same recurrence
relation, but different initial terms. That means that your new
sequence keeps many of the properties of the Fibonacci sequence. The
natural number sequence is shiftable too. Starting the natural numbers
from a different index is the same as adding a constant to the natural
number sequence. The sum of two sequences procedure ties the two
sequences by the same index in some sense. The question is, why do we
tie by the same index? Why is $a(n) + b(n)$ better than $a(n) +
b(n-1)$? If both sequences $a(n)$ and $b(n)$ are
shiftable, then $a(n) + b(n)$ might be similar to $a(n) +
b(n-1)$ and might be shiftable too. In particular, the properties
of the sum might not depend as much on how the sequences are tied
through the same index. For example, if $b(n)$ is the sequence of
natural numbers then $a(n) + b(n)$ and $a(n) + b(n-1)$ just
differ by 1. The claim is: the more shiftable your sequences, the more interesting their sum might be. The shiftability considerations correlate with my
votes for interestingness in the examples above. I find the sequence $n + Fibonacci(n)$ to be the most
interesting out of the three sequences above, the sequence $n +
prime(n)$ somewhat interesting, and the sequence $prime(n) +
Fibonacci(n)$ the least interesting.

\textbf{Product of two sequences.} The product, $p(n)$, of two
sequences $a(n)$ and $b(n)$ is defined as $p(n) =
a(n)*b(n)$. First, let us consider the product when $a =
b$. Multiplying a sequence by itself is the same as squaring the
sequence.

\textbf{Square:}

\begin{itemize}
    \item If $a(n)$ are the natural numbers, then its square is: \newline
    1, 4, 9, 16, 25, 36, 49, 64, 81, 100, 121, 144, 169, $\dots$ \newline
    A000290: The squares. 
    \item If $a(n)$ are the prime numbers, then its square is: \newline
    4, 9, 25, 49, 121, 169, 289, 361, 529, 841, 961, 1369, $\dots$ \newline
    A001248: Squares of primes.
    \item If $a(n)$ are the Fibonacci numbers, then its square is: \newline
    1, 1, 4, 9, 25, 64, 169, 441, 1156, 3025, 7921, 20736, $\dots$ \newline
    A007598: F(n)\^2, where F() = Fibonacci numbers. 
\end{itemize}

\textbf{Discussion.} These are all very interesting sequences. The first example --- the
squares --- is a very basic sequence. The second example --- the squares of primes --- has no choice but to be an exciting sequence. Namely, primes are about multiplication
properties; it is expected that you would multiply this sequence by
itself and get many interesting properties. For example, the prime squares are the numbers
that have exactly three divisors. I think the squares of the Fibonacci
numbers is the least interesting sequence out of the three. In spite
of that, by itself, the Fibonacci squares are very interesting. For
example, this sequence is a linear recurrence of order 3. It
satisfies the equation: $b(n) = 2b(n-1) + 2b(n-2) -
b(n-3)$.

\textbf{Product of two different sequences:}

\begin{itemize}
    \item If $a(n)$ are the natural numbers and $b(n)$ are the prime numbers, then $p(n)$ is: \newline
    2, 6, 15, 28, 55, 78, 119, 152, 207, 290, 341, 444, 533, $\dots$ \newline
    A033286:  n*(n-th prime).
    \item If $a(n)$ are the natural numbers and $b(n)$ are the Fibonacci numbers, then $p(n)$ is: \newline
    1, 2, 6, 12, 25, 48, 91, 168, 306, 550, 979, 1728, 3029, $\dots$ \newline
    A045925:  n*Fibonacci(n).
    \item If $a(n)$ are the prime numbers and $b(n)$ are the Fibonacci numbers, then $p(n)$ is: \newline
    2, 3, 10, 21, 55, 104, 221, 399, 782, 1595, 2759, 5328, $\dots$ \newline
    A064497:  Prime(n) * Fibonacci(n).
\end{itemize}

\textbf{Discussion.} Considerations of shiftability apply to products
too. You probably can guess my votes. Out of the three sequences above
I consider the sequence $n*Fibonacci(n)$ to be the most
interesting; the sequence $n*prime(n)$ somewhat interesting and
the sequence $prime(n)*Fibonacci(n)$ not interesting. Ironically,
the least interesting sequence I submitted myself. Why I did that is a
separate strange and sentimental story, which I might tell some other
time.

To diversify my examples, I would like to have as the third case a
more complicated and a much less famous function. Namely, in this case
$f$ is the concatenation of $x$ with the reverse of $y$.

\textbf{Concatenation of a sequence element with its reverse.} Here is a puzzle for you: look at the examples below and find what is common for all the elements of all the three sequences.

\begin{itemize}
    \item If $a(n)$ are the natural numbers, then the concatenation of $a(n)$ with the reverse of $a(n)$ is: \newline
    11, 22, 33, 44, 55, 66, 77, 88, 99, 1001, 1111, 1221, 1331, $\dots$ \newline
    Almost A056524: Palindromes with even number of digits.
    \item If $a(n)$ are the prime numbers, then the concatenation of $a(n)$ with the reverse of $a(n)$ is: \newline
    22, 33, 55, 77, 1111, 1331, 1771, 1991, 2332, 2992, 3113, $\dots$ \newline
    A067087: Concatenation of n-th prime and its reverse.
    \item If $a(n)$ are the Fibonacci numbers, then the concatenation of $a(n)$ with the reverse of $a(n)$ is: \newline
    11, 11, 22, 33, 55, 88, 1331, 2112, 3443, 5555, 8998, $\dots$ \newline
    This sequence is not in the OEIS.
\end{itemize}

\textbf{Discussion.} The answer to the puzzle: all the elements of the resulting sequences are palindromes with an even number of digits. You might also have noticed that all the elements are divisible by 11. Here is another puzzle for you: why are all the elements divisible by eleven?

Now I would like to transition from puzzles to the discussion of interestingness of these sequences. The fact that I created puzzles from these sequences might make them interesting. But if you look at my puzzles closely you can see that the puzzles are really about the first sequence out of the three. Concatenation of a number with its reverse gives you a palindrome with an even number of digits. The second sequence is the subsequence of the first sequence with prime indices. Is it interesting? I am not sure. The last sequence is not in the database, and I do not plan to submit it. You can guess why I do not want to submit it --- I really think it is not interesting.

\textbf{Concatenation of one sequence element with the reverse of another
sequence element.} Now let us go back to two variables. Suppose
$b(n)$ is different from $a(n)$. The concatenation result
depends on the order of the sequences. Obviously, $f(b(n),a(n))$
is the reverse of $f(a(n),b(n))$. For this reason I am showing
only one example out of the two for each pair of sequences:

\begin{itemize}
    \item If $a(n)$ are the natural numbers and $b(n)$ are the prime numbers, then $f(a,b)$ is: \newline
    12, 23, 35, 47, 511, 631, 771, 891, 932, 1092, 1113, $\dots$ 
    \item If $a(n)$ are the natural numbers and $b(n)$ are the Fibonacci numbers, then $f(a,b)$ is: \newline
    11, 21, 32, 43, 55, 68, 731, 812, 943, 1055, 1198, $\dots$ 
    \item If $a(n)$ are the prime numbers and $b(n)$ are the Fibonacci numbers, then $f(a,b)$ is: \newline
    21, 31, 52, 73, 115, 138, 1731, 1912, 2343, 2955, 3198, $\dots$ 
\end{itemize}

\textbf{Discussion.} The sequences above are not in the OEIS. There are
two good reasons they might not be that interesting, both of which
we have encountered before. The first reason: the prime number
and the Fibonacci number sequences are not strongly related to their
indices. The second reason: the concatenation and the reversion are
not extremely interesting operations. The main reason why they are not
interesting is that they are heavily related to the base-10
representation of numbers. In our case the sequences themselves are
not related to their base representation at all, which makes my
examples especially artificial. I have to admit that that was my goal in choosing this particular ``random'' function
--- to have very artificial examples.

\textbf{Function acting on many sequences elementwise.} Of course, as
you can guess, we can expand our definition to an integer function of
many integer variables. In this case we need many sequences to plug
in. Because I do not want to go too far away from my initial plan to start
with one or two sequences, I will give only one example here
--- the sum of my three basic sequences:

\begin{itemize}
    \item If $a(n)$ are the natural numbers, $b(n)$ are the prime numbers and $c(n)$ are the Fibonacci numbers, then $a(n) + b(n) + c(n)$ is: \newline
    4, 6, 10, 14, 21, 27, 37, 48, 66, 94, 131, 193, 287, $\dots$ 
\end{itemize}

\textbf{Discussion.} This sequence is not in the database and probably
it shouldn't be. I tried this sequence with the Superseeker and found the suggested description. 
The fact that the Superseeker
can recognize this sequence is another reason for me not to submit it.

\section{Set Operations}

In this chapter I discuss a parallel between sequences and sets. Given
a sequence, we can correspond the set of values of this sequence to
the sequence itself. Given a set of integers bounded from below, we
can create a sequence by putting the numbers in this set in increasing
order. Let us consider the set of natural numbers, which is conveniently bounded from below. That means that we can correspond a sequence to any non-empty subset of this set. And vice versa, we can correspond a subset of the set of natural numbers to any sequence of natural numbers. Note that strictly increasing sequences of natural numbers are in one-to-one correspondence with non-empty subsets of natural
numbers.

Using the described correspondence with sets we can apply set operations to sequences. In the
definitions below, I assume that $a(n)$ and $b(n)$ are
sequences of natural numbers (not necessarily increasing). To apply a set operation to sequences we first take the subsets of the natural numbers that correspond to the initial sequences, apply our set operation to them, then take the corresponding sequence as the result. Here we consider the analogs of the following set operations for sequences: complement, intersection and union.

\begin{itemize}
\item The complementary sequence $comp(n)$. Given a sequence
$a(n)$, $comp(n)$ is the sequence of natural numbers that do
not belong to $a(n)$. 
\item The intersection of two sequences $int(n)$. Given sequences
$a(n)$ and $b(n)$, the intersection $int(n)$ is the
sequence of natural numbers that belong to both $a(n)$ and
$b(n)$.
\item The union of two sequences $u(n)$. Given sequences
$a(n)$ and $b(n)$, the union $u(n)$ is the sequence of
natural numbers that belong to either $a(n)$ or $b(n)$.
\end{itemize}

Note that sometimes a set operation can produce an empty set. In this case the corresponding operation on sequences is not defined.

One of my basic sequences, the sequence of all natural numbers,
corresponds to the universal set under set operations. As a result the complement of this sequence is not defined. The union of the natural number sequence with any sequence is the natural number sequence. The intersection of the natural number sequence with a sequence $b$ is the sequence of elements $b$ put in increasing order. In particular, the intersection of the natural numbers with prime numbers is the sequence of prime numbers and the intersection of the natural number sequence with the Fibonacci sequence is a trimmed Fibonacci sequence, where we have to remove the first duplicate 1. In the examples below I omit the natural number sequence, as I just fully described its behavior under set operations.

Also, it is not very interesting to discuss the intersection or the union of a sequence with itself. The intersection or the union of a strictly increasing sequence with itself is the same sequence. In general the intersection or the union of a sequence with itself is the sequence of the elements of the original sequence in increasing order. Below I present the leftover examples of set operations applied to my basic sequences.

\textbf{Complement:}

\begin{itemize}
    \item If $a(n)$ are the prime numbers, then $comp(n)$ is: \newline
    1, 4, 6, 8, 9, 10, 12, 14, 15, 16, 18, 20, 21, 22, 24, $\dots$ \newline
    A018252: The nonprime numbers (1 together with the composite numbers of A002808).
    \item If $a(n)$ are the Fibonacci numbers, then $comp(n)$ is: \newline
    4, 6, 7, 9, 10, 11, 12, 14, 15, 16, 17, 18, 19, 20, 22, 23, $\dots$ \newline
    A001690: Non-Fibonacci numbers.
\end{itemize}

\textbf{Discussion.} The prime number sequence is property based ---
it is the sequence of all the numbers that have the property of being
prime. It is very natural to define the prime number sequence through
its corresponding set. Namely, we can define the set of prime numbers
first; then the prime number sequence is the corresponding sequence. With the
Fibonacci sequence the situation is quite opposite. The Fibonacci
sequence itself is more primary than the corresponding set. Naturally,
for property based sequences the set operations are usually more
interesting. In this case, the set of non-prime numbers is easily
defined through its property. If we exclude 1, the set of non-prime
numbers gets its own name: composite numbers. The non-Fibonacci
numbers are much less interesting.

\textbf{Intersection.} Here is my only leftover intersection example:

\begin{itemize}
    \item If $a(n)$ are the prime numbers and $b(n)$ are the Fibonacci numbers, then $int(n)$ is: \newline
    2, 3, 5, 13, 89, 233, 1597, 28657, 514229, 433494437, $\dots$ \newline
    A005478: Prime Fibonacci numbers.
\end{itemize}

\textbf{Discussion.} In general, I find the intersection operation more
interesting than the union operation. I find the intersection
especially interesting when we are dealing with property based
sequences. In this case, the intersection means numbers that have both
properties. For example, here is a very interesting intersection sequence of
numbers that are square and triangular at the same time:

\begin{itemize}
    \item If $a(n)$ are the square numbers and $b(n)$ are the triangular numbers, then $int(n)$ is: \newline
     1, 36, 1225, 41616, 1413721, 48024900, 1631432881, $\dots$ \newline
    A001110: Numbers that are both triangular and square: a(n) = 34a(n-1) - a(n-2) + 2..
\end{itemize}

\textbf{Union.} Here is my union example:

\begin{itemize}
    \item If $a(n)$ are the prime numbers and $b(n)$ are the Fibonacci numbers, then $u(n)$ is: \newline
    1, 2, 3, 5, 7, 8, 11, 13, 17, 19, 21, 23, 29, 31, 34, $\dots$ \newline
    A060634: Union of Fibonacci numbers and prime numbers.
\end{itemize}

\textbf{Discussion.} Theoretically the union is dual to the
intersection. Namely, the union is the complement of the intersection
of the complements of the given sequences. One might argue that due to
this symmetry the union should be as interesting as the
intersection. However, when we are describing the interestingness of
the sequences, very often the primary sequences are more interesting
than their complements, and the duality argument is lost. For property
based sequences the union means numbers that have either property. If
the two properties are not related to each other it is not clear to me
why the numbers with either of the properties should be joined in one
sequence. To contradict my vote for the union not being
interesting, I present an awesome example of the union of two
sequences. In this case the properties are related and the union has
dozens of interesting applications:

\begin{itemize}
    \item If $a(n)$ are the square numbers and $b(n)$ are the oblong numbers, then $u(n)$ is: \newline
    1, 2, 4, 6, 9, 12, 16, 20, 25, 30, 36, 42, 49, 56, 64, $\dots$ \newline
    A002620: Quarter-squares.
\end{itemize}

\section{Function Acting on Sets}

Suppose $f$ is a function from integers to integers. Then, given
a sequence $a$, we can define a sequence $b$: $b =
f_S(a)$ as follows: take the set of numbers corresponding
to the sequence $a$, apply the function $f$ to each number
in the set, take the resulting set (remove duplicates), then take the
sequence corresponding to the result. In other words
$f_S(a)$ is the increasing sequence of all possible
numbers that we can get when applying the function $f$ to the elements of
$a$. \textbf{Note.} This operation is defined only if applying the function $f$ to the elements of
$a$ produces a set bounded from below.

If the sequence $a$ is an increasing sequence and the function
$f$ is an increasing function, then obviously applying $f$
to the set of the elements of $a$ is the same as applying
$f$ to $a$ elementwise: $f_S(a) = f(a)$.

In the section on function acting on a sequence elementwise I had 4
different functions for my examples: adding a constant, multiplying by
a constant, the reverse square and the delta-function. Given the similarity of this operation to the function acting on a sequence elementwise, it would be
consistent to use the same 4 functions here.

\textbf{Adding a constant.} Adding a constant is an increasing function. The first two basic sequences are increasing. That means that adding a constant to the set of values of these sequences is the same as adding a constant to these sequences elementwise. The Fibonacci sequence is almost an increasing sequence. I leave it to the reader to think over the slight difference in the resulting sequences cause by adding a constant to the set of Fibonacci elements as opposed to adding a constant to the Fibonacci sequence.

\textbf{Multiplying by a constant.} For obvious reasons I do not want to multiply my sets of elements of my basic functions by negative numbers. I would happily multiply them by zero. In this case, independently of my starting sequence, my resulting sequence is a delightful sequence consisting of only one element which is 0. Multiplying our basic sequences by a positive constant gives us more diverse results than multiplying them by zero, but it is very similar to the function acting on a sequence elementwise. Namely, multiplying by a positive constant is an increasing function, and the same argument as for adding a constant applies here. That is, we saw the result of multiplying by 2 for the natural number sequence and the prime number sequence before; and with a slight change we saw the result for the Fibonacci sequence too.

\textbf{Reverse square:} Let $f(k) = Reverse(k^2)$. For example:

\begin{itemize}
    \item If $a(n)$ are the natural numbers, then $b(n)$ is: \newline
    1, 4, 9, 18, 46, 52, 61, 63, 94, 121, 144, 148, 163, $\dots$ \newline
    A074896: Squares written backwards and sorted, duplicates removed.
    \item If $a(n)$ are the prime numbers, then $b(n)$ is: \newline
    4, 9, 52, 94, 121, 148, 163, 169, 925, 961, 982, 1273, $\dots$ \newline
    Not in the OEIS.
    \item If $a(n)$ are the Fibonacci numbers, then $b(n)$ is: \newline
    1, 4, 9, 46, 52, 144, 961, 1273, 1297, 5203, 6511, $\dots$ \newline
    Not in the OEIS.
\end{itemize}

\textbf{Delta function:} Let $f(k) = \delta_1(k)$. Applying this function to the set of elements of any sequence can produce a sequence of length at most 2. Such degenerate sequences are not submitted to the database. Let us see what exactly happens to our basic sequences if we apply this function to the sets of their elements:

\begin{itemize}
    \item If $a(n)$ are the natural numbers, then $b(n)$ is: \newline
    0, 1.
    \item If $a(n)$ are the prime numbers, then $b(n)$ is: \newline
    0.
    \item If $a(n)$ are the Fibonacci numbers, then $b(n)$ is: \newline
    0, 1.
\end{itemize}

\textbf{Discussion.} I wonder what is more interesting: to apply a
function elementwise or to apply it to a set. In the first case the
order of the result is defined by the order of the underlying
sequence. In the second case the order is increasing. Which order is
better? Probably it depends on the starting sequence and the function. My example of the reverse square is not interesting in any
case, so it can't help to decide.

Suppose $f$ is a function of two variables. Then, given sequences
$a$ and $b$, we can define a sequence $c$: $c =
f_S(a,b)$ as follows: take the set of numbers
corresponding to the sequence $a$ and another set corresponding
to the sequence $b$, apply the function $f$ to each pair of
numbers from the first set and the second set, take the resulting set
(remove duplicates), then take the sequence corresponding to the
result. In other words $f_S(a,b)$ is the increasing
sequence of all possible numbers of the form $f(a(n),b(m))$.

\textbf{Sum of two sets.} Let $f(x,y) = x + y$, then
$f_S(a,b)$ is the sequence of all possible sums of the
elements from the sequence $a$ and the sequence $b$. If
$a$ is the natural number sequence and $b$ is any sequence
with the smallest element $m$, then $f_S(a,b)$ is
the sequence of natural numbers starting from $m + 1$. For this
reason in my examples I skip the cases where one of the sequences is
the natural number sequence.

\begin{itemize}
   \item If $a(n)$ and $b(n)$ are the prime numbers, then the sequence of all possible sums is: \newline
    4, 5, 6, 7, 8, 9, 10, 12, 13, 14, 15, 16, 18, 19, $\dots$ \newline
    A014091: Numbers that are the sum of 2 primes.
   \item If $a(n)$ and $b(n)$ are the Fibonacci numbers, then the sequence of all possible sums is: \newline
    2, 3, 4, 5, 6, 7, 8, 9, 10, 11, 13, 14, 15, 16, 18, $\dots$ \newline
    A059389: Sums of two nonzero Fibonacci numbers.
   \item If $a(n)$ are the prime numbers and $b(n)$ are the Fibonacci numbers, then the sequence of all possible sums is: \newline
    3, 4, 5, 6, 7, 8, 9, 10, 11, 12, 13, 14, 15, 16, 18, 19, 20, 21, 22, 23, 24, 25, 26, 27, 28, 30, $\dots$ \newline
    A132147: Numbers that can be presented as a sum of a prime number and a Fibonacci number. (0 is not considered a Fibonacci number).
\end{itemize}

\textbf{Discussion.} As I have pointed frequently out the Fibonacci
numbers are more interesting as a sequence than as a set. Therefore,
operations related to sets are usually much more interesting for the
primes than for the Fibonacci numbers. Not surprisingly, the sequence
of all possible sums of the prime numbers is the most interesting of
the three above. This sequence is related to Goldbach's conjecture that every
even integer greater than 2 can be written as the sum of two
primes. The fact that Goldbach's conjecture is one of the oldest
unsolved problems in number theory and in all of mathematics makes
this sequence especially attractive and somewhat mysterious.

\textbf{Product of two sets.} Let $f(x,y) = x*y$, then
$f_S(a,b)$ is the sequence of all possible products of
the elements from the sequence $a$ and the sequence $b$. If
$a$ is the natural number sequence and $b$ is any sequence
containing 1, then $f_S(a,b)$ is the sequence of
natural numbers. Hence, the product of the natural number sequence
with itself is the natural number sequence. Also, the product of the
natural number sequence with the Fibonacci sequence is the natural
number sequence. It is easy to see that the product of the natural
number sequence and the prime number sequence is the sequence of
natural numbers starting from 2. Here are the leftover examples:

\begin{itemize}
   \item If $a(n)$ and $b(n)$ are the prime numbers, then the sequence of all possible products is: \newline
    4, 6, 9, 10, 14, 15, 21, 22, 25, 26, 33, 34, 35, $\dots$ \newline
    A001358: Products of two primes.
   \item If $a(n)$ and $b(n)$ are the Fibonacci numbers, then the sequence of all possible products is: \newline
    1, 2, 3, 4, 5, 6, 8, 9, 10, 13, 15, 16, 21, 24, $\dots$ \newline
    A049997:  a(n) = n-th number of the form F(i)*F(j), when these Fibonacci-products are arranged in order without duplicates.
   \item If $a(n)$ are the prime numbers and $b(n)$ are the Fibonacci numbers, then the sequence of all possible products is: \newline
    2, 3, 4, 5, 6, 7, 9, 10, 11, 13, 14, 15, 16, 17, 19, $\dots$ \newline
    Almost A131511: All possible products of prime and Fibonacci numbers.
\end{itemize}

\textbf{Discussion.} And again, I find all possible products of
primes to be the most interesting sequence of the three
above. These numbers even have a name for themselves --- they are
called semiprimes.

\section{Discrete Calculus}

Given a sequence $a(n)$, the analog of the integral is the
sequence $i(a(n))$, which equals the sum of the first $n$
terms of $a(n)$. This sequence is usually called the partial sums
sequence. Similarly, the analog of the derivative is the first
difference sequence: $d(a(n)) = a(n) - a(n-1)$.

\textbf{Note.} The first term of the difference sequence is
not well defined. One of the options is to start the difference
sequence from the second term. I do not like this option because I
want all of my sequences indexed in the same way. Another option is to
assume that there is a 0 before the first term of $a(n)$, thus
artificially defining the difference for the first index. I will use
this second definition.

The integral and the derivative are complementary to each other.
The partial sums and the first difference operations are
complementary to each other in the same way. That is: $i(d(a(n))) = d(i(a(n))) =
a(n)$. \textbf{Note.} This exact equality is another good
reason to prefer the second alternative for defining the initial term
for the first difference sequence. With the first definition the
equality holds up to a constant.

\textbf{Partial sums:}

\begin{itemize}
    \item If $a(n)$ are the natural numbers, then $i(a(n))$ is: \newline
    1, 3, 6, 10, 15, 21, 28, 36, 45, 55, 66, 78, 91, 105, $\dots$ \newline
    A000217: Triangular numbers: a(n) = C(n+1,2) = n(n+1)/2 = 0+1+\linebreak2+...+n. 
    \item If $a(n)$ are the prime numbers, then $i(a(n))$ is: \newline
    2, 5, 10, 17, 28, 41, 58, 77, 100, 129, 160, 197, 238, $\dots$ \newline
    A007504: Sum of first n primes. 
    \item If $a(n)$ are the Fibonacci numbers, then $i(a(n))$ is: \newline
    1, 2, 4, 7, 12, 20, 33, 54, 88, 143, 232, 376, 609, 986, 1596, $\dots$ \newline
    Almost (shifted) A000071: Fibonacci numbers - 1. 
\end{itemize}

\textbf{First difference:}

\begin{itemize}
    \item If $a(n)$ are the natural numbers, then $d(a(n))$ is: \newline
    1, 1, 1, 1, 1, 1, 1, 1, 1, 1, 1, 1, 1, 1, 1, 1, 1, 1, $\dots$ \newline
    A000012: The simplest sequence of positive numbers: the all 1's sequence.  
    \item If $a(n)$ are the prime numbers, then $d(a(n))$ is: \newline
    2, 1, 2, 2, 4, 2, 4, 2, 4, 6, 2, 6, 4, 2, 4, 6, 6, $\dots$ \newline
    Almost A001223: Differences between consecutive primes. 
    \item If $a(n)$ are the Fibonacci numbers, then $d(a(n))$ is: \newline
    1, 1, 2, 3, 5, 8, 13, 21, 34, 55, 89, 144, 233, 377,$\dots$ \newline
    Almost (shifted) A000045: Fibonacci numbers. 
\end{itemize}

\textbf{Discussion.} The sequence of natural numbers is similar to a
polynomial of order one. It is not surprising that the partial sums
operation, which is similar to the integral, produces a sequence
corresponding to a polynomial of order two. In the same way, the first
derivative of the natural number sequence is a constant sequence (similar to
polynomials of order 0). Also, you may notice that the partial sums
as well as the first difference of the Fibonacci sequence produce the
Fibonacci sequence again. That is, the Fibonacci sequence behaves with
respect to the partial sums and the first derivative operations the
same way as the exponential function behaves with respect to the integral
and the derivative. This fact is not surprising if you remember that the
Fibonacci sequence grows similarly to the exponent of the golden ratio.

An additional natural idea is to replace the addition in the partial sums
by multiplication, thus getting partial
products. \textbf{Note.} To get the multiplicative analog of
the first difference we need to replace the subtraction operation by
division. Since the integers are not a closed set under division, I
will only supply examples for the partial products.

\textbf{Partial products:}

\begin{itemize}
    \item If $a(n)$ are the natural numbers, then the partial products sequence is: \newline
    1, 2, 6, 24, 120, 720, 5040, 40320, 362880, 3628800, 39916800, $\dots$ \newline
    A000142: Factorial numbers.
    \item If $a(n)$ are the prime numbers, then the partial products sequence is:\newline
    2, 6, 30, 210, 2310, 30030, 510510, 9699690, 223092870, 6469693230, $\dots$ \newline
    A002110: Primorial numbers. 
    \item If $a(n)$ are the Fibonacci numbers, then the partial products sequence is:\newline
    1, 1, 2, 6, 30, 240, 3120, 65520, 2227680, 122522400, 10904493600, $\dots$ \newline
    A003266: Product of first n nonzero Fibonacci numbers F(1), ..., F(n). 
\end{itemize}

\textbf{Discussion.} The prime number sequence is related to
multiplicative properties of numbers, while the Fibonacci sequence is
not. This is why I find the primorial sequence much more interesting
then the partial products of the Fibonacci numbers. Clearly, I am not
the only one who finds this sequence more interesting, as it has 
its own name.

\section{Geometric Inverse Sequence}

Suppose $a(n)$ is a positive non-decreasing sequence. Let's draw a
function graph on the $x-y$ plane corresponding to the sequence
$a(n)$. This graph consists of points $(n, a(n))$. For
consistency I would like to add a point $(0, 0)$ to the
graph, which is the same as to assume that the sequence starts with
index 0 and $a(0) = 0$. I would like to connect these points into
a piecewise linear figure looking like steps from $(0, 0)$ to
infinity. First, I add horizontal segments connecting the points
$(n-1, a(n))$ and $(n, a(n))$ for $n > 0$. Then, I
add vertical segments connecting the points $(n, a(n))$ and
$(n, a(n+1))$ for $n \ge 0$. If we symmetrically flip this
drawing with respect to the angle bisector $y = x$, we will get
another drawing that looks like steps going from $(0, 0)$ to
infinity. What is the corresponding sequence? Let us denote this new
sequence as $inv(n)$. I call this sequence the geometrical
inverse of $a$. It is easy to see that $inv(n)$ is the
maximum $m$ such that $a(m) \ n$; or equivalently, the
number of elements in the sequence $a(n)$ that are less than
$n$. Obviously, $inv(inv(a)) = a$. On the picture below you
can see the geometric inverse procedure applied to the prime number
sequence:

\begin{figure}
  \centering
  \includegraphics[scale=0.5]{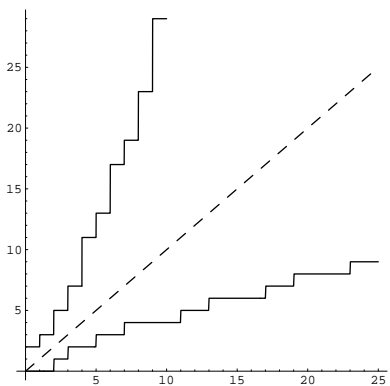}
  \caption{Geometric Inverse of Primes.}\label{basic}
\end{figure}

\textbf{Note.} The number of elements in the sequence
$a(n)$ that are less than or equal to $n$ is
$inv(n+1)$. That means that $inv(n+1)$ equals
$minimalLeftInv(n)$ defined in previous chapters for an
increasing sequence $a(n)$. Geometric inverse definition is more
general than the left inverse, as it is well defined for any
non-decreasing sequence. In particular, it is well defined for the
Fibonacci sequence.

\textbf{Geometric inverse:}

\begin{itemize}
    \item If $a(n)$ are the natural numbers, then $inv(n) = n-1$: \newline
    0, 1, 2, 3, 4, 5, 6, 7, 8, 9, 10, 11, 12, 13, 14, 15, 16, $\dots$ \newline
    A023443: n-1.
    \item If $a(n)$ are the prime numbers, then $inv(n)$ is: \newline
    0, 0, 1, 2, 2, 3, 3, 4, 4, 4, 4, 5, 5, 6, 6, 6, 6, 7, 7, $\dots$ \newline
    Almost A000720:  pi(n), the number of primes $\le$ n. 
    \item If $a(n)$ are the Fibonacci numbers, then $inv(n)$ is: \newline
    0, 2, 3, 4, 4, 5, 5, 5, 6, 6, 6, 6, 6, 7, 7, 7, 7, 7, $\dots$ \newline
    A130233:  Maximal index k of a Fibonacci number such that Fib(k)$\le$n (the 'lower' Fibonacci Inverse). 
\end{itemize}

\textbf{Discussion.} As previously mentioned, the geometric inverse of
the natural numbers and the prime numbers is the same as the
compositional minimal left inverse shifted to the right. For the
Fibonacci sequence the compositional left inverse cannot be
defined. But we presented the compositional left inverse sequences of
the trimmed Fibonacci sequence starting from the second 1. It is interesting
to compare the compositional left inverse sequences for the trimmed
Fibonacci sequence with the geometrical inverse of the Fibonacci
sequence. This comparison is left as an exercise for the reader.

The first difference of the geometric inverse shifted to the left is the
indicator sequence $ind(n)$ (also called the characteristic
function). Given a sequence $a$, the indicator sequence
$ind(n)$ equals the number of times the sequence $a$ is
equal to $n$. \textbf{Note.} For the indicator sequence,
we can remove the condition for $a(n)$ to be non-decreasing. The
necessary condition for defining the indicator function is that each
value of $a(n)$ is achieved a finite number of times.

\textbf{Indicator:}

\begin{itemize}
    \item If $a(n)$ are the natural numbers, then $ind(n) = 1$: \newline
        1, 1, 1, 1, 1, 1, 1, 1, 1, 1, 1, 1, 1, 1, 1, 1, 1, 1 $\dots$ \newline
        A000012:  The simplest sequence of positive numbers: the all 1's sequence.
    \item If $a(n)$ are the prime numbers, then $ind(n)$ is: \newline
        0, 1, 1, 0, 1, 0, 1, 0, 0, 0, 1, 0, 1, 0, 0, 0, 1, 0, $\dots$ \newline
        A010051:  Characteristic function of primes: 1 if n is prime else 0. 
    \item If $a(n)$ are the Fibonacci numbers, then $ind(n)$ is: \newline
        2, 1, 1, 0, 1, 0, 0, 1, 0, 0, 0, 0, 1, 0, 0, 0, 0, 0, 0, 0, 1, $\dots$ \newline
        A104162: Indicator sequence for the Fibonacci numbers. 
\end{itemize}

\textbf{Discussion.} You can easily prove that non-decreasing sequences are in one-to-one correspondence
with their indicators. If the sequence $a$ is strictly increasing
then its indicator takes only 0 and 1 values.

The operation of calculating the indicator function can be naturally
reversed. Here are the reverse steps: given a sequence $a(n)$,
shift it to the right, take partial sums, and then take the geometric
inverse. I call the result the reverse indicator sequence. The reverse
indicator sequence of $a(n)$ can be described as: Take $n$ 
$a(n)$ number of times.

\textbf{Reverse indicator:}

\begin{itemize}
    \item If $a(n)$ are the natural numbers, then the reverse indicator is: \newline
        1, 2, 2, 3, 3, 3, 4, 4, 4, 4, 5, 5, 5, 5, 5, 6, $\dots$ \newline
        A002024: n appears n times.
    \item If $a(n)$ are the prime numbers, then the reverse indicator is: \newline
        1, 1, 2, 2, 2, 3, 3, 3, 3, 3, 4, 4, 4, 4, 4, 4, 4, 5, $\dots$ \newline
        A083375: n appears prime(n) times. 
    \item If $a(n)$ are the Fibonacci numbers, then the reverse indicator is: \newline
        1, 2, 3, 3, 4, 4, 4, 5, 5, 5, 5, 5, 6, 6, 6, 6, $\dots$ \newline
        A072649: n occurs A000045(n) times. 
\end{itemize}

\textbf{Discussion.} I would like to draw your attention to the fact
that the reverse indicator of the Fibonacci sequence is the same
sequence as the maximal left inverse of the trimmed Fibonacci
sequence. I leave it to the reader to analyze why these sequences are
the same.

\section{Convolution of Two Sequences}

Given sequences $a(n)$ and $b(n)$ for $n$ starting with
0, their convolution is a sequence $c(n)$ defined as:
$c(0) = a(0)*b(0)$, $c(1) = a(0)*b(1) + a(1)*b(0)$, $c(2)
= a(0)*b(2) + a(1)*b(1) + a(2)*b(0)$, $\dots$, $a(n) = a(0)*b(n)
+ a(1)*b(n-1) + $\dots$ + a(n)*b(0)$, $\dots$ .

For example, if $a(n)$ is 1, 0, 0, 0, $\dots$ --- the
characteristic function of 0, then the convolution of $a(n)$ and
$b(n)$ is $b(n)$. That is, the characteristic function of 0
plays the role of the identity. For another example, if $a(n)$ is
1, 1, 1, 1, $\dots$ --- the all ones sequence, then the convolution
of $a(n)$ and $b(n)$ is the partial sums of $b(n)$
sequence. In particular, the convolution of the all ones sequence with
itself is the sequence of natural numbers shifted to the left.

Our basic sequences start with the index 1. It is easy to shift the
definition of convolution to adjust to such sequences (see Kimberling \cite{kimberling}). Given
the sequences $a(n)$ and $b(n)$ for $n$ starting with
1, the shifted convolution of them is a sequence $c(n)$ defined
as: $c(1) = a(1)*b(1)$, $c(2) = a(1)*b(2) + a(2)*b(1)$,
$c(3) = a(1)*b(3) + a(2)*b(2) + a(3)*b(1)$, $\dots$, $a(n) =
a(1)*b(n) + a(2)*b(n-1) + $\dots$ + a(n)*b(1)$, $\dots$ . In this
case the role of the identity is played by the sequence 1, 0, 0, 0,
$\dots$ --- the characteristic function of 1.

The convolution and the shifted convolution are very similar to each
other. Suppose $a(n)$ and $b(n)$ are two sequences starting
with the index 1. Suppose $a0(n)$ and $b0(n)$ are the same
sequences with 0 appended in front. Then the convolution of
$a0(n)$ and $b0(n)$ is the shifted convolution of
$a(n)$ and $b(n)$ with two zeroes appended in front. Later I
use the shifted convolution as the convolution, because our indices
start at 1.

It is easy to prove that the shifted convolution of $b(n)$ with
the natural numbers is the same as the partial sums operator applied
to the sequence $b(n)$ twice.

\textbf{Self-convolution.} Here is the shifted convolution of a basic sequence with itself:

\begin{itemize}
    \item If $a(n)$ are the natural numbers, then the self-convolution of $a(n)$ is: \newline
        1, 4, 10, 20, 35, 56, 84, 120, 165, 220, 286, 364, 455, $\dots$ \newline
        A000292: Tetrahedral (or pyramidal) numbers.
    \item If $a(n)$ are the prime numbers, then the self-convolution of $a(n)$ is: \newline
        4, 12, 29, 58, 111, 188, 305, 462, 679, 968, 1337, 1806, $\dots$ \newline
        A014342: Convolution of primes with themselves. 
    \item If $a(n)$ are the Fibonacci numbers, then the self-convolution of $a(n)$ is: \newline
        1, 2, 5, 10, 20, 38, 71, 130, 235, 420, 744, 1308, 2285, 3970, $\dots$ \newline
        A001629: Fibonacci numbers convolved with themselves. 
\end{itemize}

\textbf{Discussion.} In the OEIS database there are three natural
parameters that correlate with how interesting a sequence is:

\begin{itemize}
    \item The sequence number. Usually the more famous sequences are submitted earlier and get smaller numbers.
    \item The number of references. On each sequence page you can find a
    number in a corner in a small font that shows the number of other
    sequences referencing the given sequence. Bigger numbers usually
    correspond to more famous sequences.
    \item The size of the entry. Big entries reflect many comments and many links, and indicate an interesting sequence.
\end{itemize}

For the three sequences above, all three parameters agree. Thus, the most interesting sequence out of the three is
the sequence of tetrahedral numbers and the least interesting is the
convolution of primes with themselves.

\textbf{Convolution.} The convolution is a symmetrical operation. Here are the convolution examples for pairs of our initial sequences:

\begin{itemize}
    \item If $a(n)$ are the natural numbers and $b(n)$ are the prime numbers, then their convolution is: \newline
    2, 7, 17, 34, 62, 103, 161, 238, 338, 467, 627, 824, 1062, $\dots$ \newline
    A014148: Apply partial sum operator twice to sequence of primes.
    \item If $a(n)$ are the natural numbers and $b(n)$ are the Fibonacci numbers, then their convolution is: \newline
    1, 3, 7, 14, 26, 46, 79, 133, 221, 364, 596, 972, 1581, 2567, $\dots$ \newline
    A001924: Apply partial sum operator twice to Fibonacci numbers.
    \item If $a(n)$ are the prime numbers and $b(n)$ are the Fibonacci numbers, then their convolution is: \newline
    2, 5, 12, 24, 47, 84, 148, 251, 422, 702, 1155, 1894, 3090, $\dots$ \newline
    A023615: Convolution of Fibonacci numbers and primes.
\end{itemize}

\textbf{Discussion.} I have mentioned that the Fibonacci sequence
barely changes with respect to the partial sums operator. If we
denote the n-th Fibonacci number by $F(n)$, then the n-th partial
sum is $F(n+2) - 1$. Applying the partial sums operator
again we get a sequence whose n-th element is $F(n+4) - n
- 3$. This property is one of the reasons, why out of the
three sequences above, I find the sequence A001924
the most interesting.

\textbf{Convolutional inverse.} As I mentioned before, the sequence 1,
0, 0, 0, $\dots$ plays the role of the identity. Naturally we would wish
to define the convolutional inverse. It is easy to
see that the convolutional inverse for a sequence $a(n)$ can be
defined iff $a(1) = 1$:

\begin{itemize}
    \item If $a(n)$ are the natural numbers, then the convolutional inverse of $a(n)$ is: \newline
        1, -2, 1, 0, 0, 0, 0, 0, 0, 0, 0,  $\dots$ \newline
        Up to signs A130713: a(0)=a(2)=1, a(1)=2, a(n)=0 for n>2.
    \item If $a(n)$ are the prime numbers, then the convolutional inverse of $a(n)$ is not defined in integer sequences, but if we append primes with 1 in front, then the convolution inverse is: \newline
        1, -2, 1, -1, 2, -3, 7, -10, 13, -21, 26, -33, 53, -80, 127, $\dots$ \newline
        A030018: Coefficients in 1/P(x), where P(x) is the generating function of the primes.
    \item If $a(n)$ are the Fibonacci numbers, then the convolutional inverse of $a(n)$ is: \newline
        1, -1, -1, 0, 0, 0, 0, 0, 0, 0, 0, $\dots$ \newline
    Up to signs A130716: a(0)=a(1)=a(2)=1, a(n)=0 for n>2.
\end{itemize}

\textbf{Discussion.} The beauty of the convolution operator can be seen
if we look at the generating functions of sequences. The generating
function of the convolution of two sequences is the product of the
generating functions of these sequences. In particular, the generating
function of the convolutional inverse is the reverse of the generating
function of the sequence itself. We see that the generating functions
of the convolutional inverses of natural numbers and Fibonacci numbers
are both polynomials of order 2. Hence, the generating functions of
the natural numbers and the Fibonacci numbers are both the reverses of
second order polynomials. This means that they are both linear
recurrences of order 2. We know that fact already, but it is nice when
things come together in a different way.

\section{Binomial Transform}

There is some confusion on the web about what is called a binomial
transform. There are three different definitions very close to each
other.

Here is my first definition of a binomial transform. Given a sequence
$a(n)$ that starts with $a(0)$, the binomial transform
$b(n)$ is defined as: $b(0) = a(0)$, $b(1) = a(0) +
a(1)$, $b(2) = a(0) + 2a(1) + a(2)$, $b(3) = a(0) + 3a(1) +
3a(2) + a(3)$, $\dots$, $b(n) = a(0) + n*a(1) + $\dots$ +
C(n,k)*a(k) + $\dots$ + a(n)$, $\dots$; where $C(n,k)$ are
the binomial coefficients. This definition seems to be the most
natural out of the three. This is why it is my first choice (it is
also the first choice in Barry \cite{barry}).

\textbf{Binomial transform.} (Note that we need to add the $a(0)$ term to our initial sequences):

\begin{itemize}
    \item If $a(n)$ are the natural numbers (with 0 appended in front), then the binomial transform of $a(n)$ is: \newline
        0, 1, 4, 12, 32, 80, 192, 448, 1024, 2304, 5120, 11264, $\dots$ \newline
        A001787: n*2\^(n-1).
    \item If $a(n)$ are the prime numbers (with 1 appended in front), then the binomial transform of $a(n)$ is: \newline
        1, 3, 8, 21, 54, 137, 342, 837, 2006, 4713, 10882, 24771, $\dots$ \newline
    A030015: Binomial transform of {1, primes}. 
    \item If $a(n)$ are the Fibonacci numbers (with 0 appended in front), then the binomial transform of $a(n)$ is: \newline
        0, 1, 3, 8, 21, 55, 144, 377, 987, 2584, 6765, 17711, 46368, $\dots$ \newline
    A001906: F(2n) = bisection of Fibonacci sequence: a(n)=3a(n-1)-a(n-2). 
\end{itemize}

The reverse operation to the binomial transform is called the inverse
binomial transform. Given a sequence $a(n)$ that starts with
$a(0)$, the inverse binomial transform $b(n)$ is defined as:
$b(0) = a(0)$, $b(1) = -a(0) + a(1)$, $b(2) = a(0)
- 2a(1) + a(2)$, $b(3) = -a(0) + 3a(1) - 3a(2)
+ a(3)$, $\dots$, $b(n) = (-1)^n a(0) +
(-1)^{n-1}n*a(1) + $\dots$ + (-1)^{n-k}C(n,k)*a(k) +
$\dots$ + a(n)$, $\dots$; where $C(n,k)$ are the binomial
coefficients. \textbf{Note.} The inverse binomial transform is
called the
binomial transform at Math World \cite{BTMW}.

\textbf{Inverse binomial transform.} (Note that we need to add the $a(0)$ term to our initial sequences):

\begin{itemize}
    \item If $a(n)$ are the natural numbers (with 0 appended in front), then the inverse binomial transform of $a(n)$ is: \newline
        0, 1, 0, 0, 0, 0, 0, 0, 0, 0, 0, 0 $\dots$ \newline
        A063524: Characteristic function of 1.
    \item If $a(n)$ are the prime numbers (with 1 appended in front), then the inverse binomial transform of $a(n)$ is: \newline
        1, 1, 0, 1, -2, 5, -14, 37, -90, 205, -442, 899, -1700, 2913, $\dots$ \newline
    A030016: Inverse binomial transform of {1, primes}.
    \item If $a(n)$ are the Fibonacci numbers (with 0 appended in front), then the inverse binomial transform of $a(n)$ is: \newline
        0, 1, -1, 2, -3, 5, -8, 13, -21, 34, -55, 89, -144, 233, -377, $\dots$ \newline
    A039834: a(n+2)=-a(n+1)+a(n) (signed \hyphenation{Fi-bo-nacci} numbers); or Fibonacci numbers (A000045) extended to negative indices. 
\end{itemize}

Here is the third definition of the binomial transform. Given a
sequence $a(n)$ that starts with $a(0)$, the binomial
transform $b(n)$ is defined as: $b(0) = a(0)$, $b(1) =
a(0) - a(1)$, $b(2) = a(0) - 2a(1) + a(2)$,
$b(3) = a(0) - 3a(1) + 3a(2) - a(3)$, $\dots$,
$b(n) = a(0) - n*a(1) + $\dots$ + (-1)^k C(n,k)*a(k)
+ $\dots$ + (-1)^n a(n)$, $\dots$; where $C(n,k)$
are the binomial coefficients. This binomial transform differs only by
signs from the inverse binomial transform. That is, the n-th element of
the binomial transform by this definition is equal
$(-1)^n$ times the n-th element of the inverse binomial
transform. The beauty of this third definition is that this transform
is self-inverse. \textbf{Note.} This binomial transform is
called the
binomial transform at wiki \cite{BTwiki}.

\textbf{Binomial transform III.} (Note that we need to add $a(0)$ term to our initial sequences):

\begin{itemize}
    \item If $a(n)$ are the natural numbers (with 0 appended in front), then the third binomial transform of $a(n)$ is: \newline
        0, -1, 0, 0, 0, 0, 0, 0, 0, 0, 0, 0 $\dots$ 
    \item If $a(n)$ are the prime numbers (with 1 appended in front), then the third binomial transform of $a(n)$ is: \newline
        1, -1, 0, -1, -2, -5, -14, -37, -90, -205, -442, -899, -1700, $\dots$
    \item If $a(n)$ are the Fibonacci numbers (with 0 appended in front), then the third binomial transform of $a(n)$ is: \newline
        0, -1, -1, -2, -3, -5, -8, -13, -21, -34, -55, -89, -144, -233, $\dots$
\end{itemize}

\section{Examples: Combining Different Methods}

Here I present some examples of combining different
methods. For my examples I have chosen two sequences:

\begin{itemize}
\item Twin primes --- a very famous sequence
\item $a(n)$ is the number of n-digit powers of 2 --- not so famous a sequence, but my last OEIS submission before starting this paper.
\end{itemize}

\textbf{First example.} There are many ways to get to the twin primes from the prime number sequence. As a first step we will get to the sequence of the lesser of the twin primes:

\begin{itemize}
\item Starting with the prime numbers, take its first difference following the first definition: \newline
        1, 2, 2, 4, 2, 4, 2, 4, 6, 2, 6, 4, 2, 4, 6, 6, $\dots$ \newline
        A001223:  Differences between consecutive primes. 
\item Apply the $\delta_2$ function to every element. This is the same as taking the composition with the sequence of all zeroes except 1 in the second place: \newline
        0, 1, 1, 0, 1, 0, 1, 0, 0, 1, 0, 0, 1, 0, 0, 0, 1, 0, 0, 1, 0, $\dots$ \newline
        A100821: a(n) = 1 if prime(n) + 2 = prime(n+1), otherwise 0.
\item Apply the reverse indicator: \newline
        2, 3, 5, 7, 10, 13, 17, 20, 26, 28, 33, 35, 41, 43,  $\dots$ \newline
        A029707:  Numbers n such that the n-th and the (n+1)-st primes are twin primes.
\item Make the composition of this sequence with the prime numbers. That is, take prime numbers with the indices in the sequence above: \newline
        3, 5, 11, 17, 29, 41, 59, 71, 101, 107, 137, 149, 179,  $\dots$ \newline
        A001359: Lesser of twin primes. 
\end{itemize}

Here is another way, suggested by Alexey Radul, to get to the lesser of the twin primes sequence:

\begin{itemize}
\item Starting with the prime numbers, subtract 2 from every element: \newline
        0, 1, 3, 5, 9, 11, 15, 17, 21, 27, 29, 35, 39, 41, 45, 51, $\dots$ \newline
        A040976:  n-th prime - 2. 
\item Intersect this sequence with the prime numbers: \newline
        3, 5, 11, 17, 29, 41, 59, 71, 101, 107, 137, 149, 179,  $\dots$ \newline
        A001359: Lesser of twin primes. 
\end{itemize}

Now, from the lesser of the twin primes sequence we want to get to all
the twin primes. There are many ways to do this as well. For example:

\begin{itemize}
\item Starting with the lesser of the twin primes sequence, add 2 to every element: \newline
        5, 7, 13, 19, 31, 43, 61, 73, 103, 109, 139, 151, 181, 193, $\dots$ \newline
        A006512: Greater of twin primes.  
\item Take the union of the lesser of the twin primes with the greater of the twin primes: \newline
        3, 5, 7, 11, 13, 17, 19, 29, 31, 41, 43, 59, 61, 71, 73, 101,  $\dots$ \newline
        A001097: Twin primes. 
\end{itemize}

\textbf{Second example.} In the second example, the final sequence $a(n)$ is
the number of n-digit powers of 2. Starting with the natural numbers,
perform the following steps:

\begin{itemize}
\item Take the first difference: \newline
        1, 1, 1, 1, 1, 1, 1, 1, 1, 1, 1, 1, 1, 1, 1, 1, 1, 1, $\dots$ \newline
        A000012:  The simplest sequence of positive numbers: the all 1's sequence.  
\item Apply the binomial transform to it: \newline
        1, 2, 4, 8, 16, 32, 64, 128, 256, 512, 1024, 2048, 4096, $\dots$ \newline
        A000079: Powers of 2: a(n) = 2\^n.
\item Continue applying the binomial transform to the previous sequence 8 times: \newline
        1, 10, 100, 1000, 10000, 100000, 1000000, 10000000, 100000000,  $\dots$ \newline
        A011557: Powers of 10.
\item Take the minimal left inverse from the powers of 2 sequence: \newline
        0, 1, 2, 2, 3, 3, 3, 3, 4, 4, 4, 4, 4, 4, 4, 4, 5, 5, 5, $\dots$ \newline
        A029837: Binary order of n: log\_2(n) rounded up to next integer. 
\item Take the composition of the previous sequence with the powers of 10: \newline
        0, 4, 7, 10, 14, 17, 20, 24, 27, 30, 34, 37, 40, 44, 47, $\dots$ \newline
        A067497: Smallest power of 2 with n+1 digits (n$\ge$0). Also n such that 1 is the first digit of 2\^n.
\item Take the first difference: \newline
        4, 3, 3, 4, 3, 3, 4, 3, 3, 4, 3, 3, 4, 3, 3, 4, 3, 3, 4, 3, 3, $\dots$ \newline
        A129344: a(n) is the number of n-digit powers of 2.
\end{itemize}

\section{Generating any Sequence}

It is easy to generate any sequence from any other sequence by using
the methods described in this paper. Suppose we want to generate a
sequence $a(n)$ from the natural numbers. For my first example,
consider the function $f$ corresponding to the sequence $a$:
$f(n) = a(n)$. Then we can get the sequence $a(n)$ by
applying the function $f$ to the natural numbers.

For my second example let us apply the delta-function to the natural
numbers to get the sequence 1, 0, 0, $\dots$ . By shifting this
sequence $n$ times to the right and multiplying it by $a(n)$
we can get the sequence which has only one non-zero term and this term
is $a(n)$ at the index $n$. Then by summing all the
resulting sequences for different $n$ we get $a(n)$.

My second example requires an infinite number of steps. My first
example requires an arbitrary function. The complexity of generating
an arbitrary function is in some sense equivalent to performing an
infinite number of steps. It might be interesting to get from one
sequence to another in a finite number of operations without using
``applying a function'' procedure. Here is one way to get from
the natural numbers to the Fibonacci sequence:

\begin{itemize}
    \item Starting with the natural numbers take the convolutional inverse: \newline
    1, -2, 1, 0, 0, 0, 0, 0, 0, 0, 0,  $\dots$ 
    \item Take the partial sums of the above sequence: \newline
    1, -1, 0, 0, 0, 0, 0, 0, 0, 0, 0,  $\dots$ 
    \item Take the partial sums again: \newline
    1, 0, 0, 0, 0, 0, 0, 0, 0, 0, 0, 0, $\dots$ 
    \item Multiply the last sequence by -1 and shift two places to the right: \newline
    0, 0, -1, 0, 0, 0, 0, 0, 0, 0, 0, 0, $\dots$ 
    \item Sum the last sequence with the sequence 1, -1, 0, 0, $\dots$ in the second step: \newline
    1, -1, -1, 0, 0, 0, 0, 0, 0, 0, 0, 0, $\dots$ 
    \item Take the convolutional inverse of the last sequence: \newline
    1, 1, 2, 3, 5, 8, 13, 21, 34, 55, 89, $\dots$ 
\end{itemize}

\section{Acknowledgements}

I am thankful to Alexey Radul for criticizing my English and my writing style the first ten drafts of this paper. Alexey's help not only improved this paper tremendously, it also changed my feelings about writing in English in general. I hope it will be easier next time. I am also thankful to Jane Sherwin for checking my English in the final draft.

\end{document}